\documentclass[preprint]{elsarticle}
\usepackage{graphicx}
\usepackage{amssymb,amsmath}
\usepackage{subfig}
\newcommand{\nc}{\newcommand}
\nc{\rnc}{\renewcommand}
\nc{\x}[1]{\mbox{#1}}           % mbox avoids separations
\nc{\hs}[1]{\hspace*{#1}}
\nc{\vs}[1]{\vspace*{#1}}
\rnc{\theta}{\vartheta}
\rnc{\rho}{\varrho}
\rnc{\epsilon}{\varepsilon}
\nc{\dd}{{\mathrm{d}}}
\nc{\ii}{{\mathrm{i}}}
\nc{\ee}{{\mathrm{e}}}
\nc{\mm}[1]{\mathbf{#1}}
\nc{\mms}[1]{\boldsymbol{#1}}
\nc{\suml}[2]{\sum \limits_{#1}^{#2}}
\nc{\intl}[2]{\int \limits_{#1}^{#2}}
\rnc{\matrix}[2]{\left[\!\!\begin{array}{#1} #2\end{array}\!\!\right]}
\rnc{\vector}[1]{\matrix{c}{#1}}
\nc{\inv}{^{-1}}
\nc{\tra}{^{\mathrm T}}
\nc{\sgn}{\mathrm{sgn}}
\nc{\kn}{{k_{\mathrm n}}}
\nc{\kt}{{k_{\mathrm t}}}
\nc{\fnl}{{\tilde f}}
\nc{\xnl}{{\tilde x}}
\nc{\xc}{x_{\mathrm{c}}}
\nc{\fc}{F_{\mathrm R}}
\nc{\fst}{F_0}
\nc{\nnn}{[1]}%{{}_{\rm{n}}}
\nc{\ttt}{[2]}%{{}_{\rm{t}}}
\nc{\ntd}{N_{\mathrm{td}}}
\rnc{\Re}[1]{\operatorname{Re}\lbrace #1 \rbrace}
\nc{\prog}[1]{{\sf{#1}}}
\nc{\matlab}{\prog{Matlab}}
\nc{\ie}{i.\,e.\,}
\nc{\eg}{e.\,g.\,}
\nc{\cf}{cf.\,}
\nc{\etal}{et~al.\,}
\nc{\z}[2]{\x{\sc #1}~{\x{\cite{#2}}}}
\nc{\zo}[1]{\x{\cite{#1}}}
\nc{\fabstand}{\,}
\nc{\fp}{\fabstand.}
\nc{\fk}{\fabstand,}
\nc{\tab}[5][tbh]{\begin{table}[#1]\centering\caption{#4\label{tab:#5}}\begin{tabular}{#2}\hline #3 \\ \hline\end{tabular}\end{table}}
\newcommand{\fss}[4][tbh]{%
 \begin{figure}[#1]
 \centering
 \if@draft
  \framebox[150mm]{\raisebox{0mm}[25mm][25mm]{\texttt{#2}}}
 \else
  \includegraphics[scale=#4]{#2}
 \fi
 \caption{#3}
 \label{fig:#2}
\end{figure}
}
\newcommand{\myf}[8][tbh]{%
    \begin{figure}[#1]
        \centering
        \subfloat[#4 \label{fig:#2}]{
            \includegraphics[width=#6\textwidth]{#2}
        }
        \hspace{0.25cm}
        \subfloat[#5 \label{fig:#3}]{
            \includegraphics[width=#7\textwidth]{#3}
        }
        \caption{#8}
    \end{figure}
}
\nc{\e}[2]{\begin{equation} #1 \label {eq:#2} \end{equation}}
\nc{\ea}[2]{
\begin{eqnarray}
#1 \label {eq:#2} \end{eqnarray}}
\nc{\eal}[3][0.0ex]{
\begin{samepage}
\begin{eqnarray*}
#2
\end{eqnarray*}
\nopagebreak[4] \vs{#1} \nopagebreak[4]\vs{-2ex} \nopagebreak[4]
\begin{eqnarray}
\label {eq:#3}
\end{eqnarray}
\end{samepage}\hs{-0.35em}}
\nc{\g}[1]{{$#1$}}
\nc{\fref}[1]{{Fig.~\ref{fig:#1}}}
\nc{\frefo}[1]{{\ref{fig:#1}}}
\nc{\frefoo}[1]{{#1}}
\nc{\frefs}[1]{{Figs.~\ref{fig:#1}}}
\nc{\tref}[1]{{Tab.~\ref{tab:#1}}}
\nc{\trefo}[1]{{\ref{tab:#1}}}
\nc{\trefs}[1]{{Tab.~\ref{tab:#1}}}
\nc{\erefn}[1]{{Eq.~(\ref{#1})}}
\nc{\eref}[1]{{Eq.~(\ref{eq:#1})}} \nc{\erefo}[1]{(\ref{eq:#1})}
\nc{\erefs}[2]{{Eqs.~(\ref{eq:#1}) and (\ref{eq:#2})}}
\nc{\sref}[1]{{Section~\ref{sec:#1}}}
\nc{\srefo}[1]{\ref{sec:#1}}
\nc{\srefs}[1]{{Sections~\ref{sec:#1}}}
\nc{\aref}[1]{{{Appendix~\ref{asec:#1}}}}
\nc{\arefo}[1]{{\ref{asec:#1}}}
\nc{\arefs}[1]{{{Appendices~\ref{asec:#1}}}}
\nc{\ssref}[1]{{Subsection~\ref{sec:#1}}}
\nc{\ssrefo}[1]{\ref{sec:#1}}
\nc{\ssrefs}[1]{{Subsections~\ref{sec:#1}}}
\begin{document}

\begin{frontmatter}

\title{A HIGH-ORDER HARMONIC BALANCE METHOD FOR SYSTEMS WITH DISTINCT STATES}

\author[ids]{Malte Krack\corref{cor1}}
\ead{krack@ila.uni-stuttgart.de}
\author[ids]{Lars Panning-von Scheidt}
\author[ids]{J\"org Wallaschek}

%\address[ids]{Institute of Dynamics and Vibration Research,
%Leibniz Universit\"at Hannover, 30167 Hannover, Germany}

\cortext[cor1]{Corresponding author}

\begin{abstract}
\textit{A pure frequency domain method for the computation of
periodic solutions of nonlinear ordinary differential equations
(ODEs) is proposed in this study. The method is particularly
suitable for the analysis of systems that feature distinct states,
\ie where the ODEs involve piecewise defined functions. An
event-driven scheme is used which is based on the direct calculation
of the state transition time instants between these states. An
analytical formulation of the governing nonlinear algebraic system
of equations is developed for the case of piecewise polynomial
systems. Moreover, it is shown that derivatives of the solution of
up to second order can be calculated analytically,
making the method especially attractive for design studies.\\
The methodology is applied to several structural dynamical systems
with conservative and dissipative nonlinearities in externally
excited and autonomous configurations. Great performance and
robustness of the proposed procedure was ascertained.}
\end{abstract}


\end{frontmatter}

\section{Introduction\label{sec:introduction}}
In the fields of science and engineering, a common task is the
calculation of periodic solutions of nonlinear ordinary differential
equations. In our study, we will focus on ODEs of arbitrary
dimension involving generic, \ie possibly strong and non-smooth
nonlinear functions. In particular, we will address systems that can
comprise distinct states so that the nonlinear functions are only
piecewise defined. In mechanical engineering, such nonlinearities
arise \eg in structural systems with contact joints, where stick,
slip and lift-off are often considered as locally distinct states
\zo{john1989}. In electrical engineering, examples for such systems
are electrical circuits, where \eg transistors, rectifiers and
switches induce distinct system states. A rheological example are
superelastic shape memory alloys where the phases and phase
transformations between \eg martensite and austenite phase can be
regarded as distinct states \zo{schm2004}. Of course, many other
examples can be found in various fields of science and
engineering.\\
In order to find periodic solutions to such problems, analytical
approaches are often not applicable and computational methods have
to be employed. Besides the family of time integration methods, so
called frequency domain methods are commonly used due to their often
superior computational efficiency. The basic idea of frequency
domain methods is to choose a truncated Fourier ansatz for the
dynamic variables, thereby exploiting the periodic nature of the
solution. This class of methods gives rise to nonlinear algebraic
systems of equations. Depending on whether the solution is sought in
the frequency domain or in a collocated time domain, and whether the
residual is formulated in the frequency or time domain, the methods
can be grouped into (Multi-)Harmonic Balance Method
\zo{urab1965,nayf1979}, Trigonometric Collocation Method
\zo{sall2011} and Time Spectral Method \zo{gopi2005}. Among these
methods, the Multi- or High-order Harmonic Balance Method (HBM) is
probably the most commonly applied method.\\
For the HBM, it is generally necessary to compute the spectrum of
the nonlinear function that governs the ODE. This task can generally
be performed by different methods. In the following, we will focus
on those methods that are capable of treating systems with distinct
states.\\
The Alternating-Frequency-Time (AFT) scheme \zo{came1989} is one of
the most commonly applied approaches in this context. The AFT scheme
involves a sampling of the nonlinear function and subsequent
back-transformation into frequency domain. Advantages of this method
are the broad applicability, the comparatively small implementation
effort and the low computational effort for evaluating the residual
function. The latter aspect is particularly true if the (Inverse)
Fast Fourier Transform is used for the transformation between time
and frequency domain. A drawback is that nonlinearities with
distinct states involve special treatment. A sampling of the
nonlinear function is not straight-forward, because the current
state at a specific time instant is not always a priori known. So
called predictor-corrector schemes \zo{guil1998} are frequently
employed to perform the switching between different states for these
systems. In classical AFT schemes, the sampling points are fixed,
and do not need to coincide with the state transition time instants.
This inherently induces discretization errors. Hence, the
sensitivity of the transition time instants with respect to
arbitrary parameters cannot be captured accurately, resulting in
inaccurate derivatives, in particular for higher order
derivatives.\\
More recently, a purely frequency-based formulation was proposed by
\textsc{Cochelin and Vergez} \zo{coch2009}. The authors applied the
Asymptotic Numerical Method to expand the periodic solution into a
power series based on high-order derivatives of the nonlinear
function. In order to obtain these derivatives efficiently, a so
called quadratic recast is performed where the original system of
equations is transformed into a system of only quadratic order. An
advantage of this method is the computationally robust and efficient
continuation of the solution. A drawback is obviously the required
quadratic recast which can be difficult for generic types of
nonlinear functions. Moreover, systems with distinct states need to
be artificially smoothed in order to accomplish a closed-form
quadratic recast. This smoothing procedure induces
inaccuracies compared to the original non-smooth model.\\
In order to avoid the shortcomings of a required recast or the
degenerated accuracy due to sampling, a pure frequency domain
formulation for the original system with distinct states can instead
be used. Such an approach necessarily involves the direct
calculation of the transition time instants between the states. For
high-order HBM, these approaches have only been developed for
special types of nonlinearities so far. For example \z{Petrov and
Ewins}{petr2003b} developed an analytical formulation of the HBM for
piecewise linear friction interface elements in structural dynamical
problems. In this study, the approach in \zo{petr2003b} is extended
to generic systems with an arbitrary number of distinct states, see
\sref{methods_of_analysis}. Analytical formulations can be developed
in case of piecewise polynomial systems, as it will be shown in
\ssref{pp}. Moreover, the formulation facilitates the analytical
calculation of gradients of up to second order as an inexpensive
postprocessing step, see \ssref{sensana}. To demonstrate the
capabilities and the performance of the proposed methodology,
several numerical examples are studied in \sref{numerical_examples}.
Finally, conclusions are drawn in \sref{conclusions}.

\section{Methods of Analysis\label{sec:methods_of_analysis}}
%This section is organized as follows: In \ssref{hbm}, the Harmonic
%Balance Method is applied to a system with distinct states. The
%formulation involves the calculation of the state transition time
%instants, which is addressed in \ssref{statetrans}. The proposed
%method is then applied to the class of piecewise polynomial systems
%in \ssref{pp}. Finally, analytical expressions for the first and second order
%sensitivities of the solution are derived in \ssref{sensana}.

\subsection{Harmonic Balance Method for systems with distinct states\label{sec:hbm}}
Consider a system whose dynamics can be described by a first-order
ordinary differential equation,
\e{\dot{\mm y} = \mm f\left(\mm y,t\right)\fk}{nonlinear_ode}
in which \g{\dot{()}} denotes derivative with respect to time $t$.
It is assumed that the generally nonlinear function \g{\mm f} is
piecewise defined within closed regions of the state space of \g{\mm
y}. These closed regions in state space are denoted \textit{states}
throughout this paper. These states shall not be confused with the
vector \g{\mm y} which is sometimes also referred to as state in
literature since it represents a point in state space.
\fss[thb]{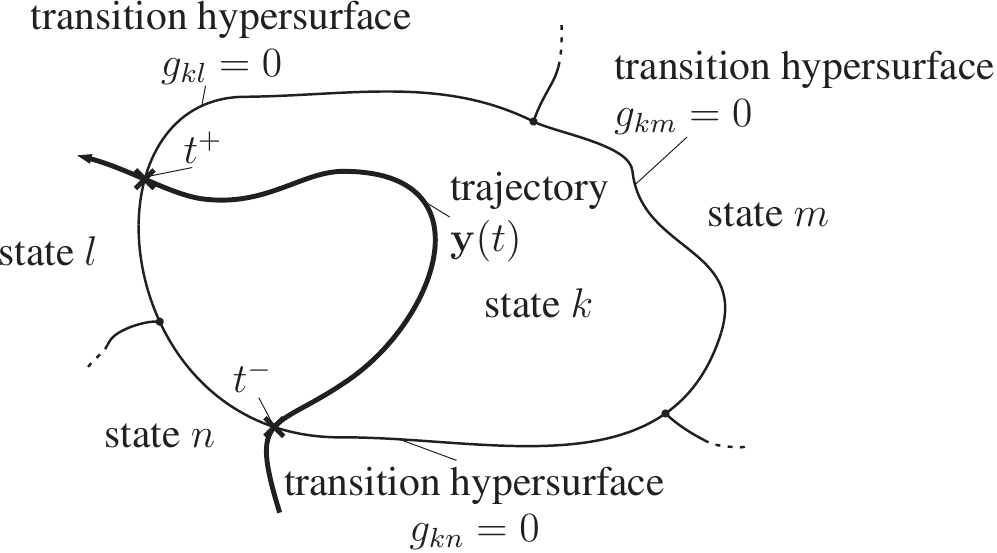}{Illustration of the dynamics of a system with
distinct states}{0.75}
\\
As time evolves, the system can assume several states, see
\fref{fig1}. A transition between these states is termed state
transition and the corresponding time instant is called state
transition time instant in the following. The system enters a
specific state at the corresponding transition time \g{t^-} and
leaves it at \g{t^+}. Each possible state $k$ consists of a
nonlinear function \g{\mm f_k}, transition conditions \g{g_{kl}}
which roots define a transition hypersurface to state $l$, and
internal variables \g{\mm v_k}:
\ea{\nonumber \text{\underline{Definition of `state $k$'}} && \\
\nonumber \text{Nonlinear function:} & \mm f_k\left(\mm y(t),\mm
v_k,t\right)\fk &\\
\nonumber \text{Transition conditions:} & g_{kl}\left(\mm y(t),\mm
v_k,t\right)\,, & \forall l\in\mathcal L_k\fk\\
\text{Internal variables:} & \mm v_k\left(\mm y(t^-),\mm
f(t^-)\right)\fp}{state}
The set \g{\mathcal L_k} is a set of integers indicating a possible
next state, the system can assume after being in state $k$.\\
It is assumed that the function \g{\mm f} is smooth within a state
and continuous at the state transitions. It should be noted at this
point that the advantages of the proposed method can be particularly
exploited for the case of piecewise polynomial systems, as it will
be shown in \ssref{pp}, although the derivations in the following
are not restricted to these.\\
Internal variables are introduced in the state definition
\erefo{state} to facilitate the treatment of hysteresis effects. In
a hysteretic system, the dynamics do not explicitly depend on the
current value \g{\mm y(t)} but on the time history of \g{\mm y}.
Internal variables can therefore be used to carry this
history-dependent effect over the state hypersurface, which
manifests itself in the dependence of the nonlinear function \g{\mm
f_k} and the transition hypersurface \g{g_{kl}} on \g{\mm v_k}, see
\eref{state}. Note that hysteretic systems will also be addressed in
the numerical examples. For systems without these effects, of
course, the
introduction of internal variables is not necessary.\\
Periodic, steady-state solutions to \eref{nonlinear_ode} are sought
in this study. To this end, the High-order Harmonic Balance Method
can be applied \zo{nayf1979}. Hence, a Fourier series truncated to
harmonic order \g{H} represents the ansatz for the dynamic variables
\g{\mm y(t)},
\e{\mm y(t) \approx \suml{n=-H}{H}{\mm Y_n \ee^{\ii n\Omega
t}}}{fourier_ansatz}
Herein, \g{\Omega} is the fundamental angular frequency of the
response and \g{\ii = \sqrt{-1}} is the imaginary unit. The Fourier
coefficients \g{\mm Y_n} are symmetric, \g{\mm Y_{-n} =
\overline{\mm Y}_{n}}, where \g{\overline{()}} denotes complex
conjugate, since \g{\mm y(t)} is a real-valued function in time.
Substitution of \eref{fourier_ansatz} into the differential equation
\erefo{nonlinear_ode} and Fourier-Galerkin projection \zo{urab1965}
gives rise to a nonlinear algebraic system of equations in the
unknowns \g{{\mm Y}_{n}} and \g{\tau_j^-,\tau_j^+},
\ea{\text{solve} & \ii n\Omega\mm Y_n - \mm F_n\left(\mm Y_{-H},
\cdots, \mm Y_{H}\right) = \mm 0\,,\quad
n=-H,\cdots,H\fk\label{eqm_hbm}\\
\text{with} & \mm F_n = \frac{1}{2\pi}\intl{(2\pi)}{}\mm f(\mm
y,\tau)\ee^{-\ii
n\tau}\dd\tau=\frac{1}{2\pi}\suml{j=1}{J}\intl{\tau_j^-}{\tau_j^+}\mm
f(\mm y,\tau)\ee^{-\ii n\tau}\dd\tau\fk \label{fnl_hbm}\\
\text{subject to} &
\tau_J^+=\tau_1^-+2\pi\,,\quad\tau_j^+=\tau_{j+1}^-\,\,\,\,
\forall\,\, j=1\cdots J\fp\label{pstates_hbm}}{hbm}
For convenience, the normalized time \g{\tau=\Omega t} has been
introduced. During one period of oscillation, the system assumes a
total number of $J$ states. It should be emphasized that neither the
set of states nor the state transition time instants
\g{\tau_j^-,\tau_j^+} are a priori known. As indicated in the
constraint \erefn{pstates_hbm}, the \g{\tau_j^-,\tau_j^+} are
periodic and continuous to cover an entire time period, as a
consequence of the periodic ansatz. In this study, it is proposed to
directly compute the periodic set of transition time instants
\g{\tau_j^-,\tau_j^+}, which is developed in the following
subsection.\\
Once the transition time instants are known for given \g{\mm Y_n},
the integrals in \erefn{fnl_hbm} can be evaluated to formulate the
residual in \erefn{eqm_hbm}. Owing to the piecewise definition of
the function \g{\mm f}, it is convenient to split up the integral in
\erefn{fnl_hbm} into $J$ summands, where each of the summands is an
integral with the transition time instants as integral limits.

\subsection{Periodic set of state transition time instants\label{sec:statetrans}}
\fss[h!]{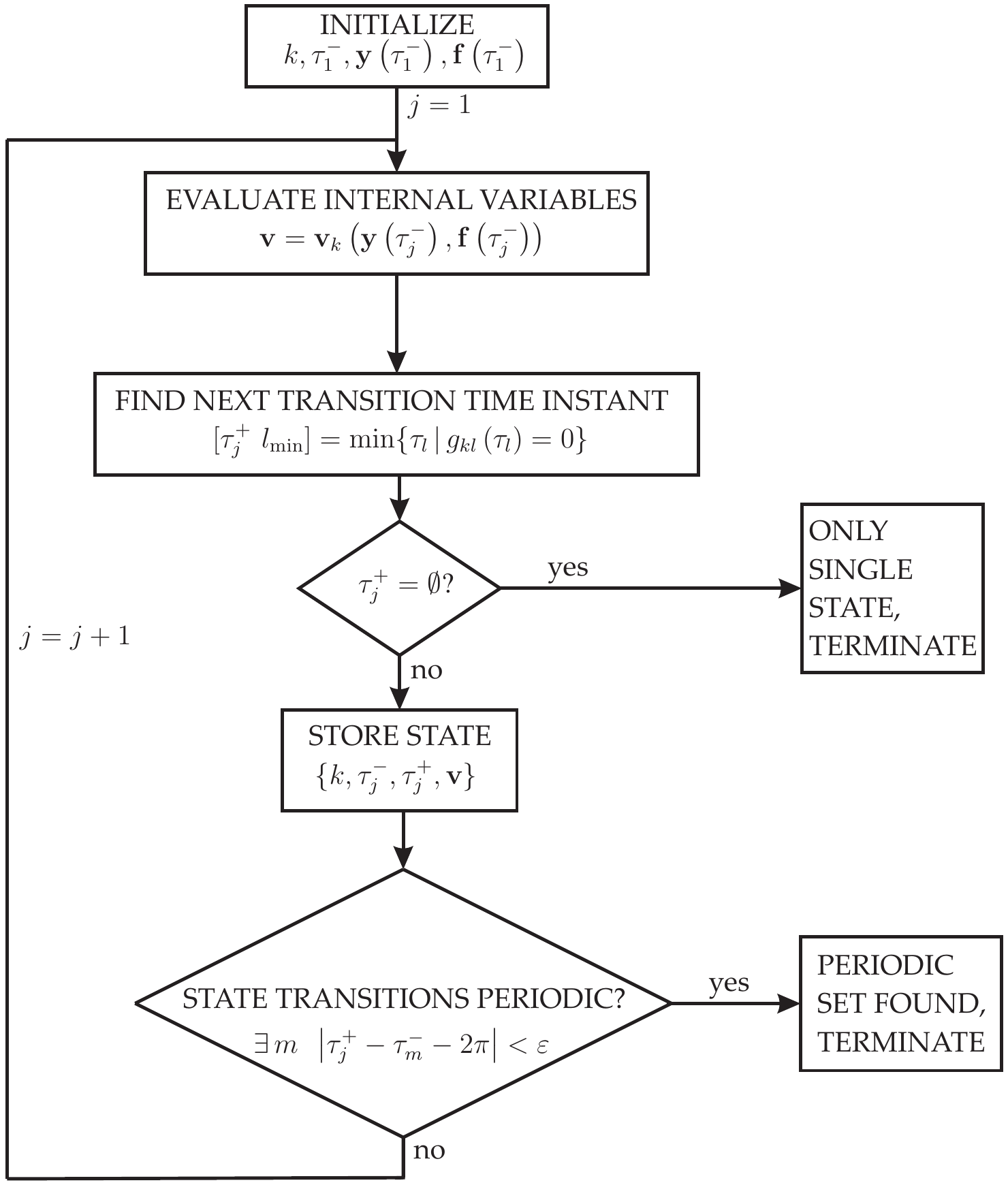}{Algorithm for the calculation of the periodic set of
state transition time instants}{0.75}
In \fref{fig2}, an algorithm is summarized that is capable of
finding a periodic set of state transition time instants for an
arbitrary system with or without distinct states. Starting from an
initial time \g{\tau_1^-}, state $k$ and according initial function
value \g{\mm f(\tau_1^-)}, the next states are iteratively computed
until a periodic set of state transitions is found. It is therefore
assumed that a periodic set of state transitions exists and the
algorithm is attracted to it. During the numerical studies, no case
was observed where this assumption was disproved.\\
After evaluating the internal variables, the next roots \g{\tau_l}
of all possible state transition conditions \g{g_{kl}} are computed
and the minimum is taken. In the special case when there is no next
state, the system remains in this state for all times and the
algorithm terminates. Note that this includes the special class of
systems with only a single state.\\
If a next state exists, the current state is stored for subsequent
evaluation of the Fourier coefficients. If the state transitions are
periodic - according to a specified tolerance $\varepsilon$ - the
algorithm can terminate, otherwise $j$ is incremented and the loop
is repeated.
%        - when the state at a particular time instant can be
%        determined and this state does not contain any internal
%        variables

\subsection{Computation and continuation of the solution}
In general, the solution to
\x{Eqs.~(\ref{eqm_hbm})-(\ref{pstates_hbm})} cannot be obtained in
closed form and an iterative numerical procedure has to be employed
instead. In this study, a Newton-Raphson method combined with a
predictor-corrector continuation scheme was used \zo{seyd1994}. The
numerical performance of the solution procedure was greatly enhanced
by providing analytically calculated gradients of the residual, as
derived in the following subsection.

\subsection{Analytical calculation of gradients and sensitivities of
the solution \label{sec:sensana}} Gradients of the residual are
often required in a numerical solution procedure for the algebraic
system of equations in \erefn{eqm_hbm}. Moreover, higher-order
derivatives at the solution point can be used to expand the solution
in a Taylor series. An approximate solution thus becomes available
in the vicinity of the current solution point in parameter space
without the need for re-computation. The Taylor expansion with
respect to the unknown variables can be employed as a predictor in a
numerical continuation procedure. Taylor expansions with respect to
system parameters are particularly interesting for
parametric studies, uncertainty analysis and optimization.\\
In this study, the analytical calculation of gradients of first and
second order is presented. We focus on the Fourier coefficients of
the nonlinear function \g{\mm f} since the sensitivities of the
other term in \erefn{eqm_hbm} is straight-forward. The first and
second order sensitivities of \g{\mm F_n} read
\ea{\frac{\partial\mm F_n}{\partial\psi} =
\suml{j=1}{J}\intl{\tau_j^-}{\tau_j^+}\frac{\partial \mm
f}{\partial\psi}\ee^{-\ii n\tau}\dd\tau + \mm
f(\tau_j^+)\frac{\partial\tau_j^+}{\partial\psi} - \mm
f(\tau_j^-)\frac{\partial\tau_j^-}{\partial\psi}=
\suml{j=1}{J}\intl{\tau_j^-}{\tau_j^+}\frac{\partial \mm
f}{\partial\psi}\ee^{-\ii n\tau}\dd\tau\fk\label{eq:testbild}\\
\frac{\partial^2\mm F_n}{\partial\phi\partial\psi} =
\suml{j=1}{J}\intl{\tau_j^-}{\tau_j^+}\frac{\partial^2 \mm
f}{\partial\phi\partial\psi}\ee^{-\ii n\tau}\dd\tau +
\frac{\partial\mm
f(\tau_j^+)}{\partial\psi}\frac{\partial\tau_j^+}{\partial\phi} -
\frac{\partial\mm
f(\tau_j^-)}{\partial\psi}\frac{\partial\tau_j^-}{\partial\phi}\fp\quad}{fnl_sensana}
Herein, \g{\phi,\psi} are arbitrary scalar variables such as the
components of the Fourier coefficients \g{\mm Y_n}, the frequency
\g{\Omega} or any system parameter. The Leibniz integral rule was
applied to derive \erefs{testbild}{fnl_sensana} since the integral
limits might and often do depend on the parameters. Note that the
last two summands in the first-order sensitivity cancel each other
out in the sum over one period due to the assumed continuity of
\g{\mm f} and
the periodicity condition in \erefn{pstates_hbm}.\\
The calculation of the derivative of the function \g{\mm f} within a
state is typically straight-forward. In contrast, the sensitivities
of a transition time instant \g{\tau_j} is more complex. It has to
be calculated by implicit differentiation of the active transition
condition \g{g(\tau_j)=0}. The resulting first- and second-order
sensitivities of the transition time instants read
\ea{\nonumber \frac{\partial\tau_j}{\partial\psi} &=&
\left[\frac{\partial g}{\partial\tau}\right]^{-1}\frac{\partial g}{\partial\psi}\fk\\
\frac{\partial^2\tau_j}{\partial\phi\partial\psi} &=&
\left[\frac{\partial
g}{\partial\tau}\right]^{-1}\left[\frac{\partial^2
g}{\partial\phi\partial\psi} + \frac{\partial
g}{\partial\phi\partial}\frac{\partial\tau_j}{\partial\psi} +
\frac{\partial
g}{\partial\psi\partial}\frac{\partial\tau_j}{\partial\phi}
\frac{\partial^2
g}{\partial\tau^2}\frac{\partial\tau_j}{\partial\phi}
\frac{\partial\tau_j}{\partial\psi}\right]\fp}{tau_sensana}
In \eref{tau_sensana}, all functions are evaluated at the transition
time instant \g{\tau_j}. It should be remarked that the time
derivative of the transition condition \g{\frac{\partial
g}{\partial\tau}} is nonzero at a regular zero crossing so that the
inverse in \eref{tau_sensana} is well-defined. Note that only
first-order derivatives of the state transition time instants
\g{\tau_j} are directly included in \eref{fnl_sensana}. However,
second-order derivatives may be required for the calculation of the
sensitivities of the internal variables \g{\mm
v_k\left(\mathbf{y}\left(\tau_j^-\right),
\mathbf{f}\left(\tau_j^-\right)\right)}, see definition
\erefo{state}.

\subsection{Application to piecewise polynomial systems\label{sec:pp}}
All previous developments are valid for the class of piecewise
smooth systems. In the sequel of this study, we will focus on the
large subclass of piecewise polynomial systems. For this class, all
functions \g{\mm f_k, g_{kl}, \mm v_k} are polynomials in the
components of \g{\mm y}, which makes the efficient formulation of
the previously derived expressions particularly cheap. In order to
solve \x{Eqs.~(\ref{eqm_hbm})-(\ref{pstates_hbm})} the basic
operations (a) \textit{add/subtract}, (b) \textit{multiply}, (c)
\textit{integrate} and (d) \textit{calculate roots} are required to
find the periodic set of state transitions and to carry out the
integration indicated in
\erefn{fnl_hbm}. These operations can be directly performed in Fourier space.\\
The \textit{multiplication} of two scalar functions
\g{a(\tau),b(\tau)} with associated Fourier coefficients \g{\mm
A=[A_{-H},\cdots,A_H],\mm B=[B_{-H},\cdots,B_H]} can be expressed as
a convolution in Fourier space,
\e{\mathcal F\lbrace a\cdot b\rbrace = \mm A \ast\mm
B\fp}{fourier_product}
Herein, \g{\mathcal F} indicates the Fourier Transform and \g{\ast}
denotes convolution. Note that powers of a Fourier series can be
calculated by recursive multiplication.\\
The \textit{integration} of a truncated Fourier series in the time
interval \g{\tau^-} to \g{\tau^+} can be expressed as follows (see
\arefo{fourier_integral}):
\e{\intl{\tau^-}{\tau^+} a(\tau)\ee^{-\ii n\tau}\dd\tau =
\left(\tau^+-\tau^-\right) A_n + \suml{m=-H,m\neq
n}{H}\frac{\ee^{\ii\left(m-n\right)\tau^+}-
\ee^{\ii\left(m-n\right)\tau^-}}{\ii\left(m-n\right)}A_m\fp}{fourier_integral}
This equation can be applied to the evaluation of the integrals in
\erefn{fnl_hbm}.\\
There are efficient as well as robust numerical methods for the
\textit{calculation of the roots} of a truncated Fourier series, see
\eg \zo{boyd2006}. Most of these methods simply compute the roots of
the associated complex polynomial in \g{z=\ee^{\ii\tau}}. Such
methods are available in many computational software frameworks like
\matlab.\\
It is important to note that the harmonic order is increased by the
convolution in \eref{fourier_product}, \ie when products or powers
of a Fourier series are generated. It is therefore proposed to
truncate the Fourier series of the nonlinear function \g{\mm F_n} to
the original order $H$ in \erefn{eqm_hbm} so that the resulting
number of equations is equal to the number of unknowns.

\subsection{On the numerical performance and accuracy of the proposed method\label{sec:numperformance}}
In all numerical studies of the piecewise polynomial systems
presented in \sref{numerical_examples}, the computational bottleneck
was observed to be the root finding of the complex polynomials
involved in the calculation of the state transition time instants.
State-of-the art polynomial root finding algorithms are based on the
computation of the eigenvalues of a so called companion matrix, for
which the computational complexity increases approximately with the
number of harmonics cubed \g{\mathcal O(H^3)}. The interested reader
is referred to \zo{boyd2006} for a detailed analysis of the
computational cost for this operation.\\
In contrast to the root finding operation, carrying out time-domain
integration, differentiation and multiplication by means of
summation and matrix multiplication according to the derived
closed-form expressions in \ssref{pp} is comparatively efficient.
This has some noteworthy implications for the analytical
calculations of the gradients: The evaluation of first and second
order derivatives represents an efficient post-processing step,
since their calculation only involves comparatively cheap vector and
matrix multiplications, as indicated in \ssref{sensana}.\\
It should be remarked that the accuracy of the proposed method,
particularly regarding the gradients, relies on the direct
calculation of the state transitions. The conventional AFT scheme is
characterized by an inherent discretization error. This causes a
severe limitation for the achievable accuracy. In this context, it
is interesting to note that in a piecewise linear system, the
second-order derivatives essentially result from the sensitivities
of the transition time instants, which can be easily verified from
\eref{tau_sensana}. As these sensitivities are not captured by the
AFT scheme, the second-order sensitivities would be identical to
zero in this case. This emphasizes the superiority of the proposed
method with respect to the AFT scheme regarding accurate sensitivity
analysis.

\section{Numerical examples\label{sec:numerical_examples}}
We have implemented the methodology proposed in
\sref{methods_of_analysis} in a computational framework in the
\matlab software environment. We used an object-oriented software
architecture to exploit operator overloading capabilities. For
example, we defined a Fourier series class that implements the
required operations add/subtract, multiply, integrate and compute
roots. Moreover, we developed and used an Automatic Differentiation
class similar to the one described in \zo{fort2006,neid2010} to
carry out the analytical sensitivity analysis up to second order. A
database of state formulations was created that
includes the nonlinearities presented in this section.\\
The numerical examples in this study comprise structural dynamical
systems. Application of the proposed methodology to fields other
than structural dynamics, \eg electrical networks, is considered
straight-forward but beyond the scope of this study. For structural
dynamical systems, the vectors \g{\mm y,\mm f} can be written as
follows:
\ea{\mm y = \vector{\mm x\\ \mm{\dot x}}\,,\quad \mm f =
\vector{\mm{\dot x}\\ -\mm M\inv\left(\mm D\mm{\dot x}+\mm K\mm x +
\mm f_{\mathrm e}(t)+\mm{\fnl}(\mm x,\mm{\dot
x})\right)}\fp}{dynsys}
Herein, \g{\mm M,\mm D, \mm K} are structural mass, damping and
stiffness matrices, \g{\mm x} is the vector of generalized
displacements, \g{\mm f_{\mathrm e},\mm{\tilde f}} are generalized
excitation and nonlinear forces.\\
The numerical examples can be categorized in two groups. In the
examples in \ssrefs{2dof_cubic}-\ssrefo{fgl} a 2-Degree-of-freedom
(DOF) system with an attached single nonlinear element \g{\fnl} is
considered, see \fref{fig3}. The example for \ssref{beam} is a
cantilevered beam with contact constraints and will be described
later. It should be noted that the example systems with a small
number of DOFs were considered for clarity reasons. The methodology
proposed in this paper can generally be applied to systems with
arbitrary number of DOFs, including large-scale Finite Element
Models.
\fss[t]{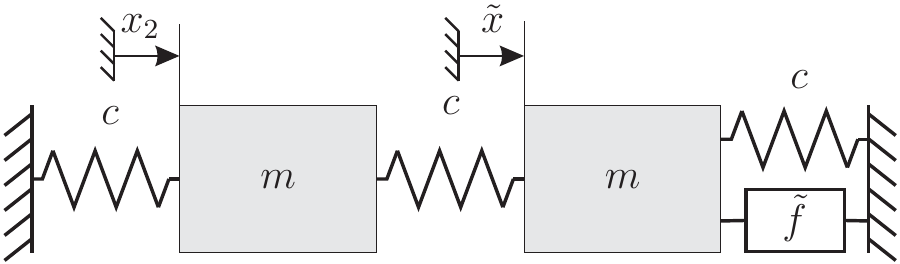}{2-DOF system with nonlinear element}{1.0}
\\
For the 2-DOF system, the structural matrices and the nonlinear
force vector have the following form:
\ea{\nonumber \mm M = \matrix{cc}{1 & 0\\ 0 & 1}\,,\,\,\mm D=\mm
0\,,\,\, \mm K = \matrix{cc}{2 & -1\\ -1 &
2}\,,\,\,\mm{\fnl}=\vector{\fnl\\ 0}\fp}{2dof}
The corresponding vector of generalized coordinates is \g{\mm x\tra
= \matrix{cc}{\xnl & x_2}}, where \g{\xnl} denotes the nonlinear
DOF.

\subsection{2-DOF system with cubic spring\label{sec:2dof_cubic}}
For a first demonstration of the methodology, a cubic spring
nonlinearity is considered. The nonlinearity can be described by a
single state without transition conditions and no internal
variables. In the notation introduced in \eref{state}, the state
definition reads as listed in \tref{2dof}.
\tab[h]{lc}{ & state 1\\
$\fnl$ & $0.5 \xnl^3$\\
$g$ & (-)\\
$v$ & (-)}{State definition of a system with cubic spring}{2dof}
\\
This example emphasizes once again that purely polynomial, \ie
smooth nonlinearities are a special case of the piecewise polynomial
class treated in this study.
\myf[t]{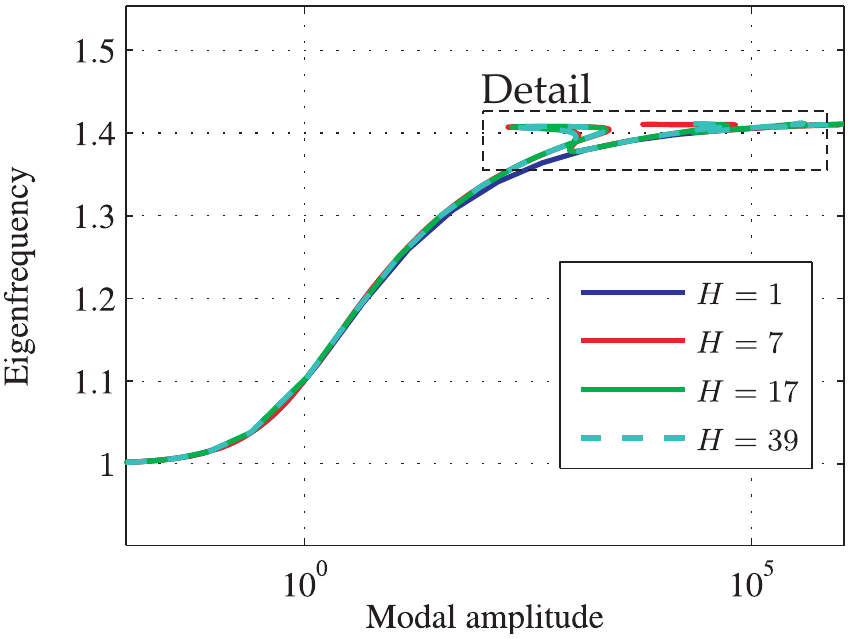}{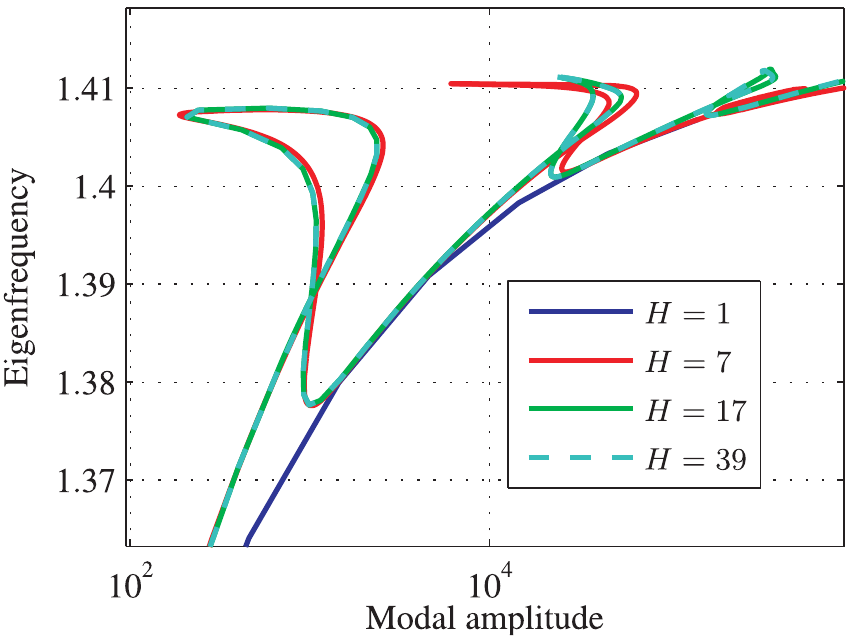} {}{}{.45}{.45}{Frequency-Energy-Plot of a
2-DOF system with cubic spring (~(a) overview, (b) detail~)}
\\
The 2-DOF system is investigated in autonomous configuration. The
proposed method was used for the calculation of the nonlinear normal
modes. Great convergence behavior was ascertained.\\
A thorough study of the nonlinear normal modes of this system can be
found in \z{Kerschen \etal}{kers2009} and shall not be repeated
here. Instead, only the so called Frequency-Energy-Plot (FEP) of the
first nonlinear mode is depicted in \frefs{fig4a}-\frefo{fig4b}.
Throughout this study, amplitude and frequency axes in the figures
are scaled by their values for the linear case. The eigenfrequency
increases with the modal amplitude due to the stiffening effect of
the cubic spring. For large amplitudes, the energy localizes in the
left mass in \fref{fig3}. The system exhibits several internal
resonances in the considered modal amplitude range \zo{kers2009},
causing so called tongues in the FEP, see \fref{fig4b}. Apparently
several harmonics are required to accurately predict the nonlinear
modal interactions.

\subsection{2-DOF system with piecewise polynomial spring\label{sec:snap_through}}
Again, the 2-DOF system is considered, however, the cubic spring is
now replaced by a piecewise polynomial spring. The
force-displacement characteristic is given by the function depicted
in \fref{fig5a}. The nonlinearity was defined by introducing three
states, each with a polynomial force \g{\fnl} as listed in
\tref{snap_through}. As it was shown in \sref{pp}, the High-order
Harmonic Balance residual equations can be formulated analytically
for this class of systems by the new technique proposed in this
paper.
\tab[h]{lccc}{ & state 1 & state 2 & state 3\\
$\fnl$ & $-(\xnl-1)^2+1$ & $-(\xnl-2)$ & $(\xnl-3)^2-1$\\
$g$ & $g_{12}=\xnl-1$ & $g_{21}=g_{12}$, $g_{23}=g_{32}$ & $g_{32}=\xnl-3$\\
$v$ & (-) & (-) & (-)}{State definition of a system with piecewise
polynomial spring}{snap_through}
\myf[b!]{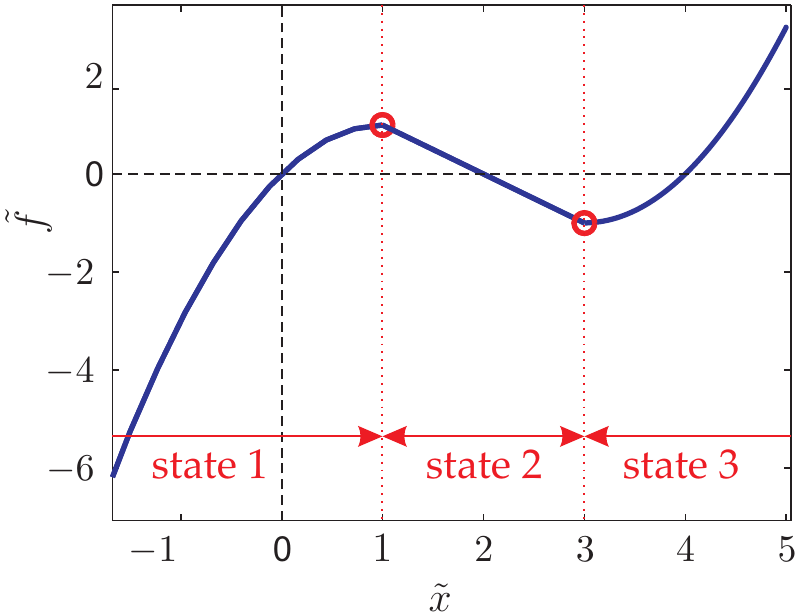}{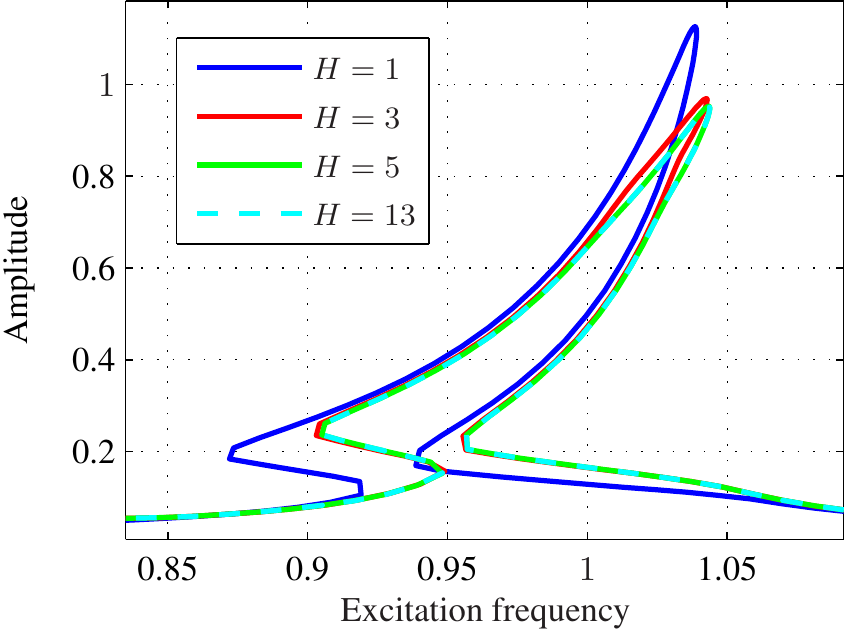}{}{}{.45}{.45}{Characteristics of a 2-DOF
system with piecewise polynomial spring (~(a) force-displacement
relationship, (b) forced response function~)}
\\
The central part of the characteristic is linear with negative
slope. The two neighboring states have a quadratic
force-displacement characteristic. Note that the piecewise
polynomial spring is conservative with a unique force-displacement
relationship. Hence, the state formulation does
not require any internal variables.\\
A harmonic force excitation at the linear mass in a frequency range
close to the first eigenfrequency is imposed. The forced response
function was calculating using the proposed method and is
illustrated in \fref{fig5b}. Overhanging branches occur in the
forced response characteristic: For moderate vibration amplitudes,
\ie for small vibrations around the equilibrium point \g{\xnl=0},
the system exhibits softening behavior and the amplitude-frequency
curve is bent to the left. The effective stiffness decreases with
increasing amplitude due to negative slope in the force-displacement
characteristic. For larger vibration amplitudes, the effective
stiffness increases due to the quadratic branches, resulting in a
stiffening behavior and the amplitude-frequency curve is bent to the
right. Apparently, several harmonics have to be considered in the
harmonic expansion to accurately predict the dynamic behavior of the
system.

\subsection{2-DOF system with elastic Coulomb friction element\label{sec:twodoffric}}
An elastic Coulomb friction or Masing element \zo{masi1923} is
attached to the 2-DOF system in \fref{fig3}. Tangential stiffness
$\kt$ and friction force limit $\fc$ characterize this nonlinearity.
The nonlinearity can assume two states: Stick (state 1) and slip
(state 2), see \tref{elastic_coulomb}.
\myf[b!]{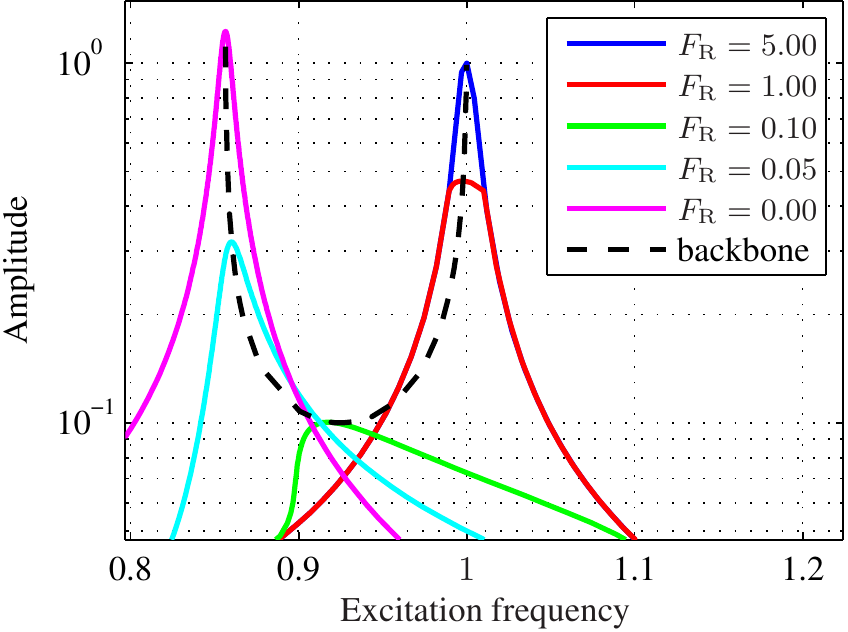}{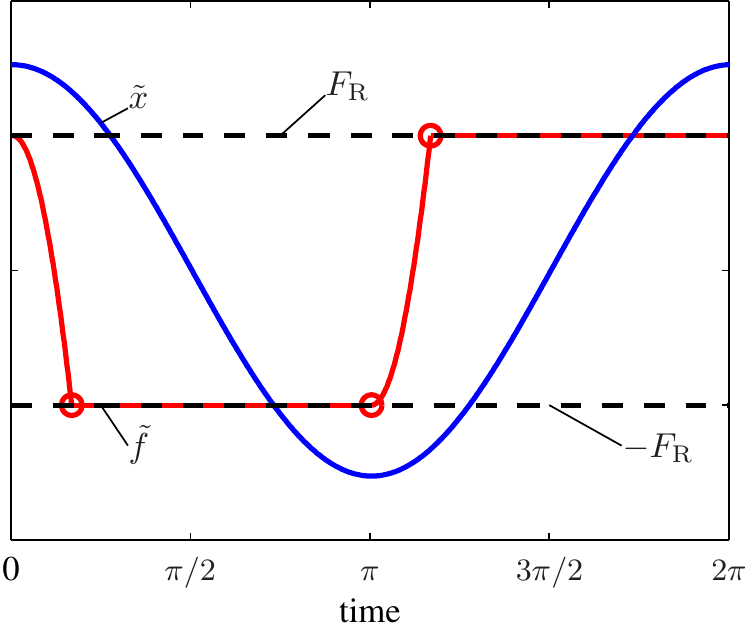}{}{}{.48}{.42}{Forced response of a system
with elastic Coulomb nonlinearity (~(a) forced response functions
for different values of the normal load, (b) typical time history~)}
\myf[t!]{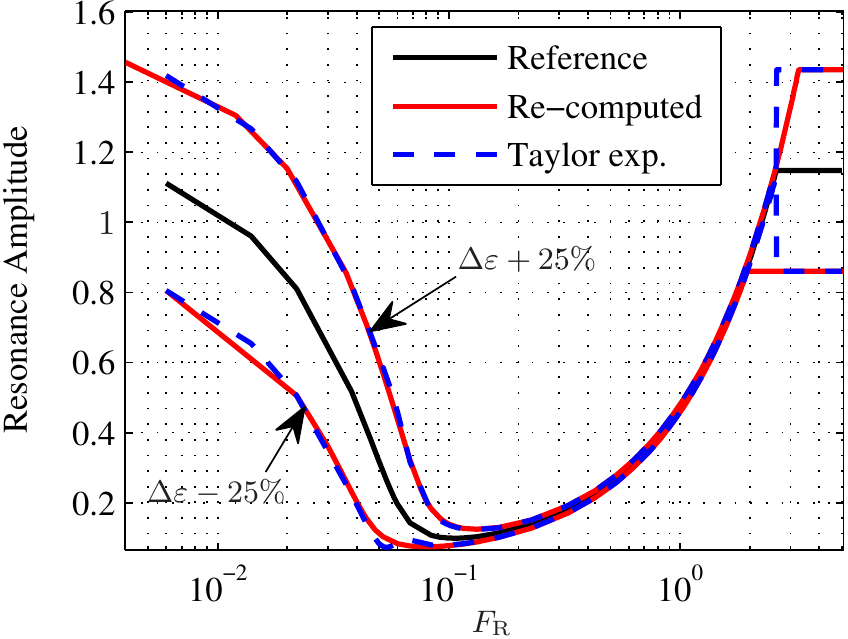}{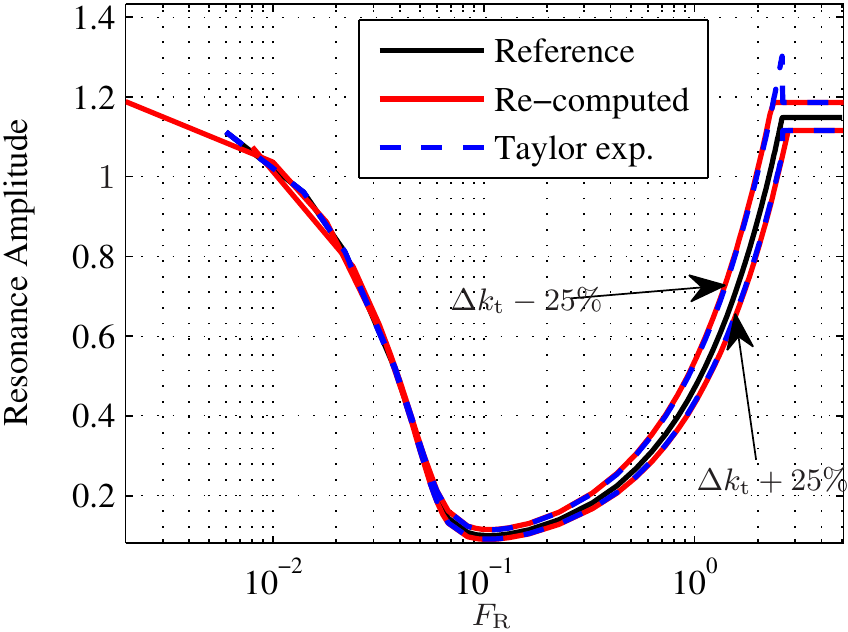}{}{}{.45}{.45}{Resonance amplitude as a
function of the normal load (~(a) variation of the excitation level
\g{\epsilon}, (b) variation of the tangential stiffness \g{\kt}~)}
\tab[h]{lcc}{ & state 1 & state 2\\
$\fnl$ & $\kt(\xnl-v_1)$ & $v_2\fc$\\
$g$ & $g_{12}=\fnl^2-(\fc)^2$ & $g_{21}=\dot\xnl$\\
$v$ & $v_1=\xnl(\tau_j^-)-\frac{\fnl(\tau_j^-)}{\kt}$ &
$v_2=\sgn\fnl(\tau_j^-)$}{State definition of a system with elastic
Coulomb friction element}{elastic_coulomb}
\\
If the elastic friction force reaches its limit value \g{\fc}, \ie
\g{g_{12}=0}, a stick-to-slip transition occurs. If a reversal point
is reached (\g{\dot{\xnl}=g_{21}=0}), a stick phase is initiated.
Internal variables for the elastic Coulomb element are the Coulomb
slider position $v_1$ and the slip direction $v_2$.
%        - note that an alternative would be to distinguish slip
%        state into positive and negative slip state
\\
Again, a harmonic force excitation is imposed at the linear mass. In
\fref{fig6a}, the forced response function in the vicinity of the
eigenfrequency of the first mode is depicted for different values of
the friction force limit \g{\fc}. A tangential stiffness value of
\g{\kt=0.35} was specified. For large values of \g{\fc}, the Coulomb
element is fully stuck so that the hysteresis degenerates to a line
and no damping effect is introduced by the friction element. For
vanishing values of \g{\fc}, the slider can slip freely so that the
hysteresis is flat and again there is no friction damping effect. In
between these extreme cases, a significant amplitude reduction due
to friction damping can be achieved. Moreover, the resonance
frequency increases as the value of \g{\fc} increases due to the
coupling effect of the friction element. The backbone curve that
connects the maxima of the forced response functions was directly
calculated by applying the strategy described in \zo{krac2012a} to
the methodology proposed in this paper. A typical time history of
both displacement \g{\tilde x} and nonlinear force \g{\tilde f} is
illustrated in \fref{fig6b}. Owing to the moderate value of \g{\kt},
the response remains essentially harmonic. The transitions between
stick and slip state can be well-observed from the time
history of the force in \fref{fig6b}.\\
The suitability of the analytically formulated sensitivities is now
investigated. To this end, the resonance amplitude of the first mode
is depicted as a direct function of \g{\fc} in
\frefs{fig7a}-\frefo{fig7b}. These so called optimization curves are
often used for design purposes, see \eg
\zo{bert1998,petr2006b,krac2012b}. In addition to the nominal
parameter set, the optimization curve is also illustrated for
slightly smaller and larger ($\pm 25\%$) excitation level and
tangential stiffness values. The results were obtained by
second-order Taylor expansion about the reference solution (Taylor
exp.). For comparison, the optimization curves were also computed
directly at the new parameter point (Re-computed). The results agree
well in a wide range of the \g{\fc} value. However, the Taylor
expansion fails in predicting the fully stuck configuration, \ie for
very high \g{\fc} values. As it was also reported in
\zo{brau1993,krac2012a}, it is not possible to accurately predict
the dynamic behavior beyond the point where the order or number of
states change.

\subsection{2-DOF system with superelastic shape memory alloy\label{sec:fgl}}
The hysteresis effect of superelastic shape memory alloys (SMA) can
be employed for damping of mechanical structures. A sophisticated
modeling approach would involve constitutive as well as
thermodynamical aspects, see \eg \zo{bern2003}. This is, however,
regarded as beyond the scope of this study and a simplified
rheological piecewise linear model \zo{schm2004} shall be considered
instead. The associated hysteresis can be described by five distinct
states as illustrated in \fref{fig8} and listed in \tref{fgl}. The
system features a purely elastic state (1). The forward and reverse
transformation between austenite and martensite phase is described
by the states (2) and (4). Beyond a certain displacement, a linear
onset (3) is used to describe the superelastic behavior. Depending
on the displacement evolution in time, an intermediate state (5) can
also be reached. Note that the point symmetry of the hysteresis is
exploited in the state definition in \tref{fgl}.
\fss[b!]{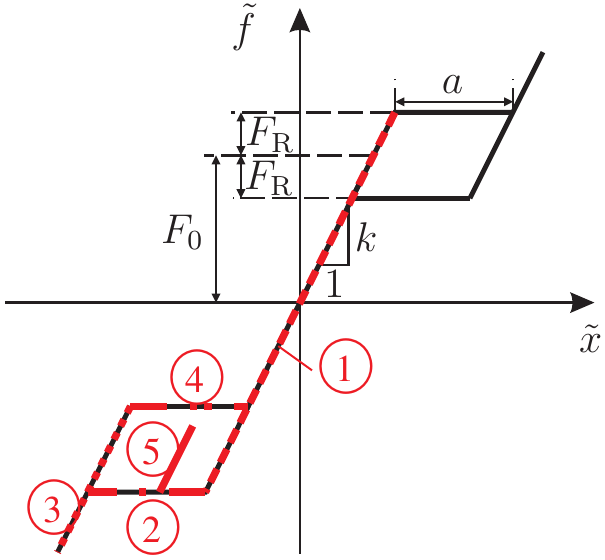}{Approximated hysteresis of a superelastic shape
memory alloy}{1.0}
\myf[t!]{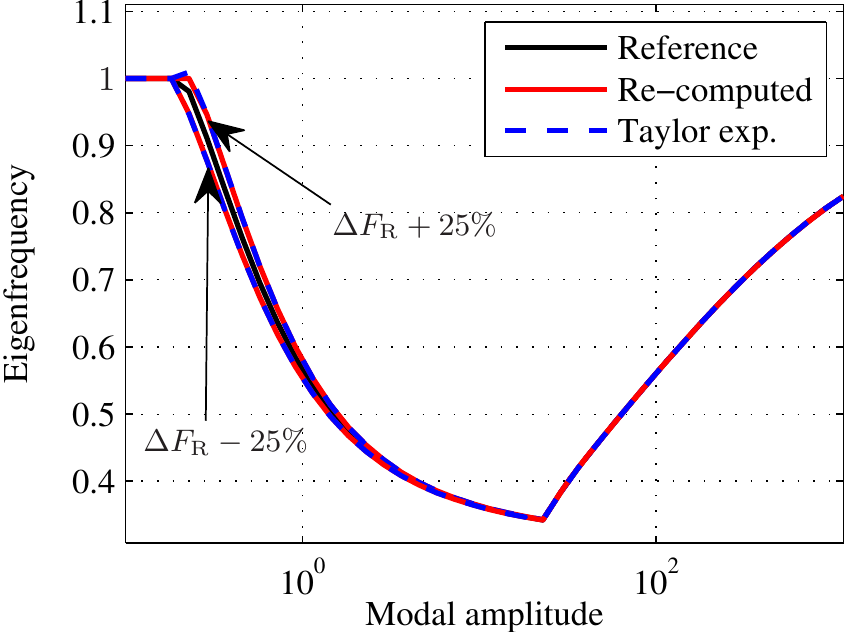}{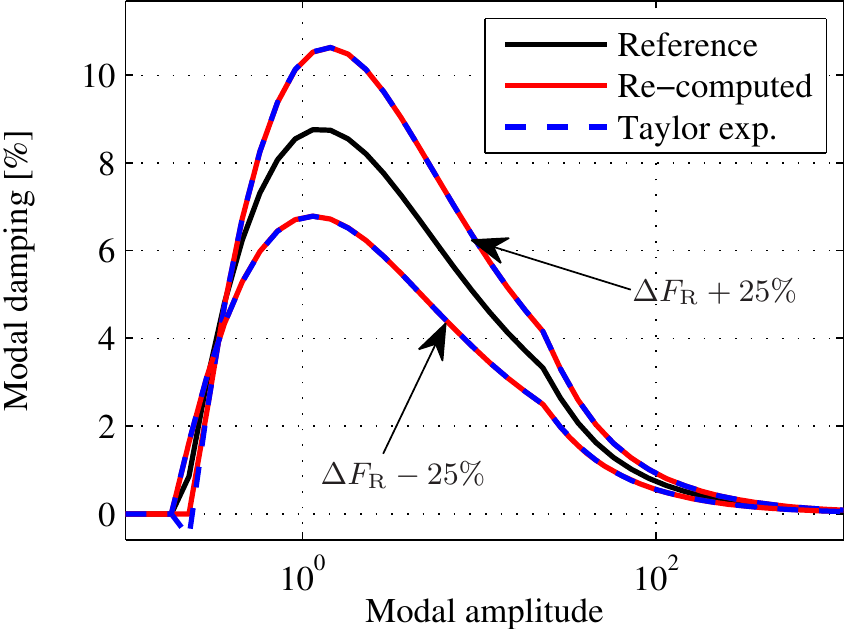}{}{}{.45}{.45}{Modal properties of a 2-DOF
system with superelastic shape memory alloy (~(a) eigenfrequency,
(b) modal damping~)}
\tab[h!]{lccc}{ & state 1 & state 2 & state 3\\
$\fnl$ & $k\xnl$ & $v_2\left(\fst+\fc\right)$ &
$k(\xnl-v_3 a)$\\
$g$ & $g_{12}=\fnl^2-(\fst+\fc)^2$ & $g_{23}=kv_2\xnl-ka-\fst-\fc$,
& $g_{34} = \fnl-v_3(\fst-\fc)$\\
& & $g_{25}=\dot\xnl$ &\\
$v$ & (-) & $v_2=\sgn\,\xnl(\tau_j^-)$ & $v_3=v_2$\\
\hline\\
\hline
& state 4 & state 5\\
$\fnl$ & $v_4\left(\fst-\fc\right)$ & $k(\xnl-v_5)$\\
$g$ & $g_{41} = kv_3\xnl-\fst+\fc$, &
$g_{52}=\fnl^2-(\fst+\fc)^2$,\\
& $g_{45} = \dot\xnl$ & $g_{54}=\fnl^2-(\fst-\fc)^2$\\
$v$ & $v_4=v_2$ & $v_5 =
\xnl(\tau_j^-)-\frac{\fnl(\tau_j^-)}{\kt}$}{State definition of a
system with superelastic shape memory alloy}{fgl}
\\
The SMA-type nonlinearity was also applied to the 2-DOF system in
\fref{fig3}. The Nonlinear Modal Analysis technique proposed in
\zo{laxa2009,krac2013d} was employed in conjunction with the
formulations of the nonlinearities proposed in this study.
Eigenfrequency and the modal damping ratio were computed with
respect to the modal amplitude of the first mode. The results are
depicted in \frefs{fig9a}-\frefo{fig9b}. For small vibration
amplitudes, the system always remains in state 1, \ie the damping
vanishes and the eigenfrequency is constant. For moderate vibration
amplitudes, the phase transformation occurs to a certain extent so
that the damping value increases and the eigenfrequency is reduced
due to the softening effect. For large vibration amplitudes, the
effect of the hysteresis cycles becomes smaller again so that
eigenfrequency and
damping value asymptotically approach their linearized values again.\\
As in the previous example, the sensitivities of the nonlinear
dynamic analysis results have been computed to formulate a
second-order Taylor series in the system parameters. Using the
sensitivity results, the modal properties have been expanded with
respect to the parameter \g{\fc}, \cf \frefs{fig9a}-\frefo{fig9b},
for \g{\pm 25\%} deviation from the nominal value. The results agree
well with the re-computed results.

\subsection{Beam with friction and unilateral contact\label{sec:beam}}
\fss[b]{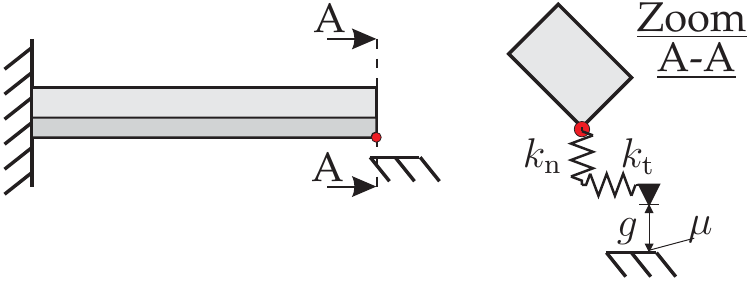}{Cantilevered beam with friction and unilateral
contact at its free end}{0.7}
As a final example, a clamped beam with combined friction and
unilateral contact was investigated. The system is depicted in
\fref{fig10}. A finite element code was used to mesh the geometry
and derive the structural matrices of the cantilevered beam for the
initial configuration. The finite element model comprised $10,098$
DOFs. A single node-to-ground contact element was attached to the
free end as depicted in \fref{fig10}. In contrast to the example in
\ssref{twodoffric}, the contact model additionally accounts for the
variation of the normal load and possible lift-off. The contact
nonlinearity can thus assume three distinct states as listed in
\tref{nt_elastic}: Separation (state 1), stick (state 2), slip
(state 3). System parameters are tangential stiffness $\kt$,
friction coefficient \g{\mu}, normal stiffness $\kn$ and normal gap
$g$.
\myf[t]{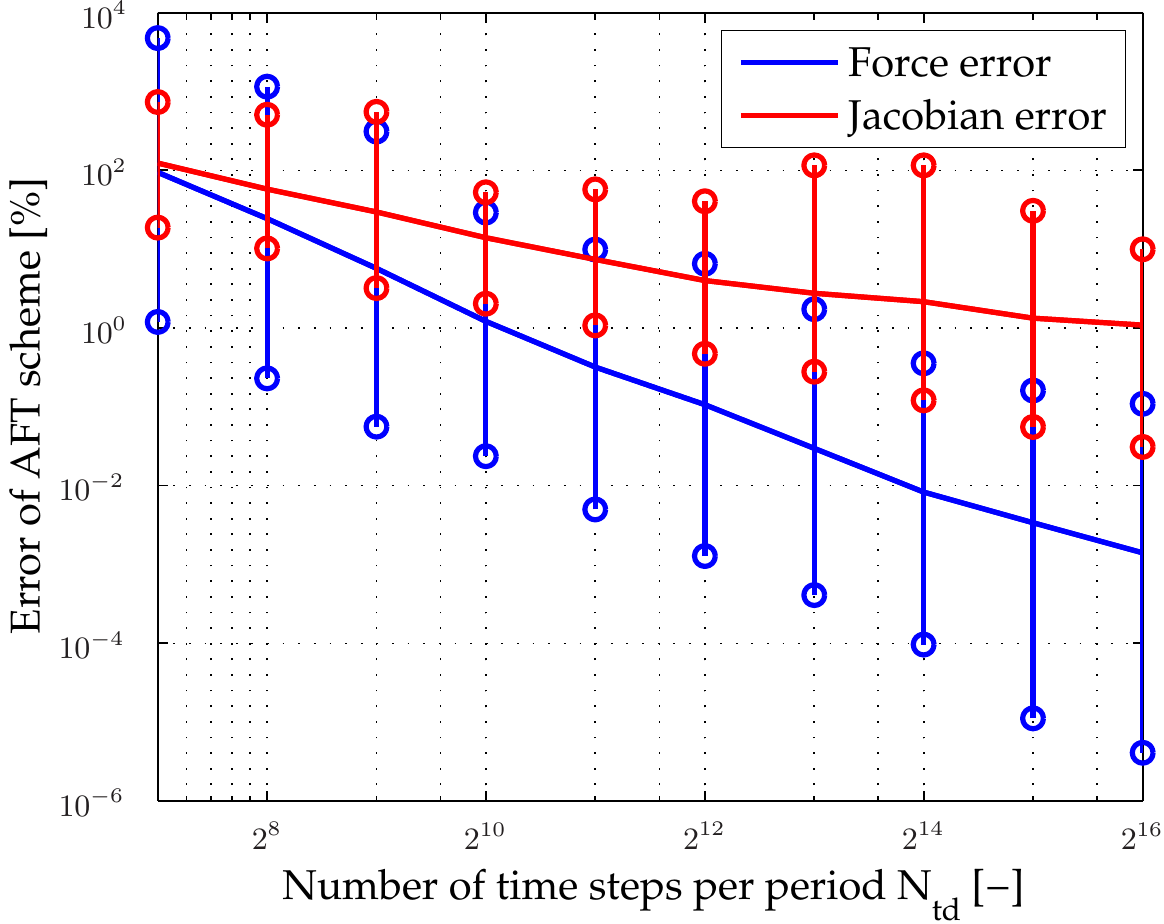}{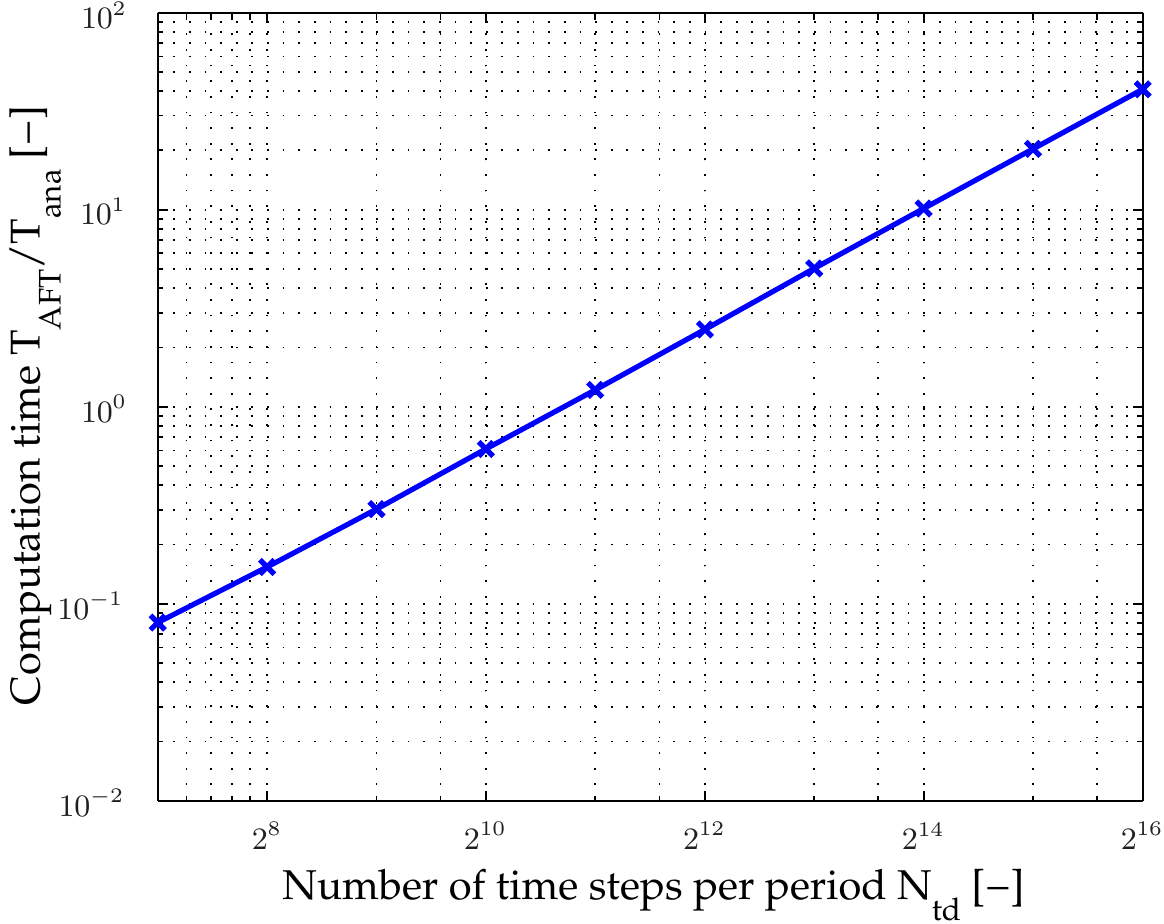}{}{}{.45}{.45}{Comparison of
Alternating-Frequency-Time scheme with proposed method (~(a)
accuracy of function and derivative, (b) computational effort~)}
\tab[h]{lccc}{ & state 1 & state 2 & state 3\\
$\mm \fnl$ &
$\vector{0\\ 0}$ & $\vector{\kn\left(\xnl\nnn+g\right)\\
\kt\left(\xnl\ttt-v_2\right)}$ & $\vector{\kn\left(\xnl\nnn+g\right)\\
v_3\mu\kn\left(\xnl\nnn+g\right)}$ \\
$g$ &
$g_{12}=\xnl\nnn+g\,\text{if}\,\kn{\dot\xnl}\nnn>\kt{\dot\xnl}\ttt$,
& $g_{21}=\xnl\nnn+g$, &
$g_{31} = g_{21}$,\\
&
$g_{13}=\xnl\nnn+g\,\text{if}\,\kn{\dot\xnl}\nnn\leq\kt{\dot\xnl}\ttt$
& $g_{23}=\left(\fnl[2]\right)^2-\left(\mu\fnl[1]\right)^2$ &
$g_{32}=\dot{\fnl}[2]-\kt{\dot\xnl}\ttt$\\
$v$ & (-) & $v_2=\xnl\ttt(\tau_j^-) - \frac{\fnl[2](\tau_j^-)}{\kt}$
& $v_3=\sgn\,\fnl[2](\tau_j^-)$}{State definition of a system with
friction and unilateral contact}{nt_elastic}
\\
For this particular example, a comparison with the conventional AFT
scheme was performed in terms of accuracy and computational effort,
\cf \frefs{fig11a}-\frefo{fig11b}. Only the nonlinear force
calculation is considered for the comparison. Random vectors of
complex displacement amplitudes was generated. A number of $1,000$
random vectors was large enough to obtain convergence of the
performance statistics. Seven harmonics have been considered in the
analysis. In \fref{fig11a}, the mean, minimum and maximum error of
the force and the Jacobian are depicted with respect to the number
of time samples $\ntd$ per period used in the AFT scheme. The
accuracy of the AFT scheme can be increased by increasing the number
of time steps. An larger number of time steps yields a better
accuracy, but also a higher computational effort. The computational
effort $T_{AFT}$ essentially increases linearly with the number of
time steps. The effort quickly exceeds the one required for the
proposed method ($T_{ana}$), \cf \fref{fig11b}. It has to be
remarked that only the nonlinear force is considered in the
performance comparison. It is expected that the resulting error in
the predicted response is less significant than that of the force or
the Jacobian.
\myf[t]{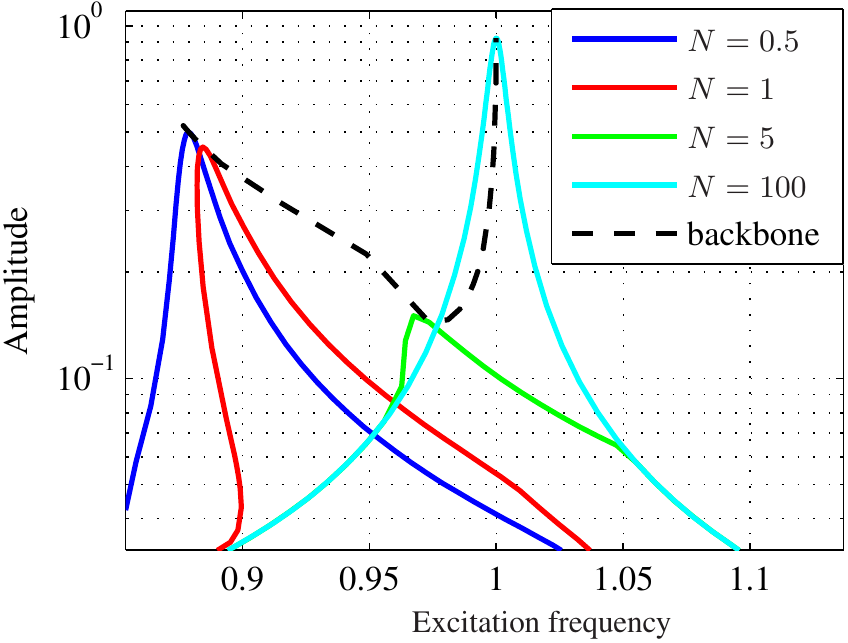}{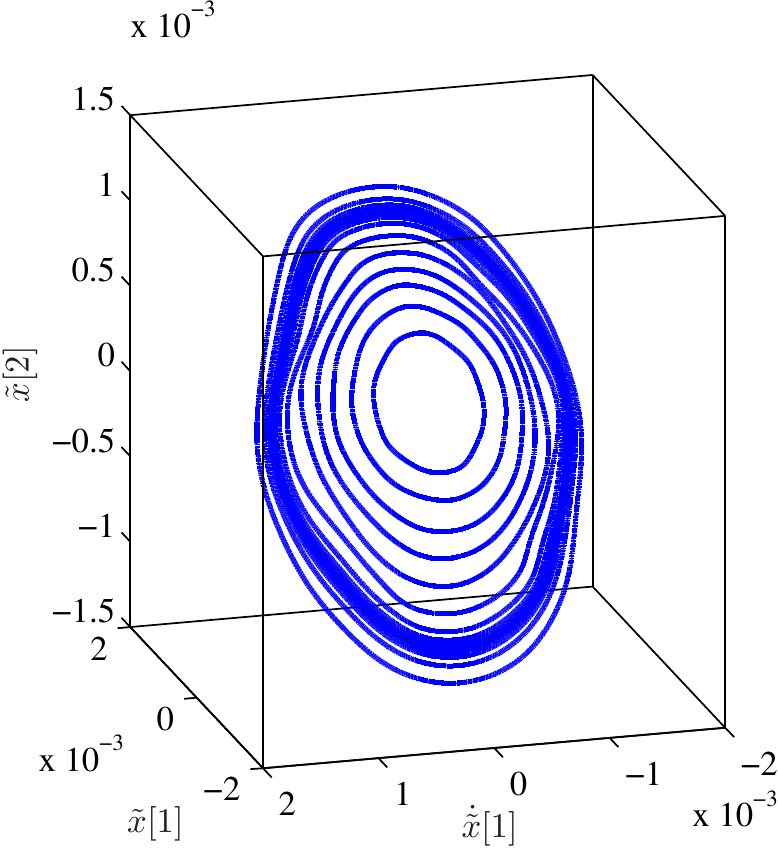}{}{}{.48}{.42}{Forced response of a beam with
friction and unilateral contact (~(a) forced response functions for
different values of the normal load, (b) orbits along the backbone
curve~)}
\myf[h!]{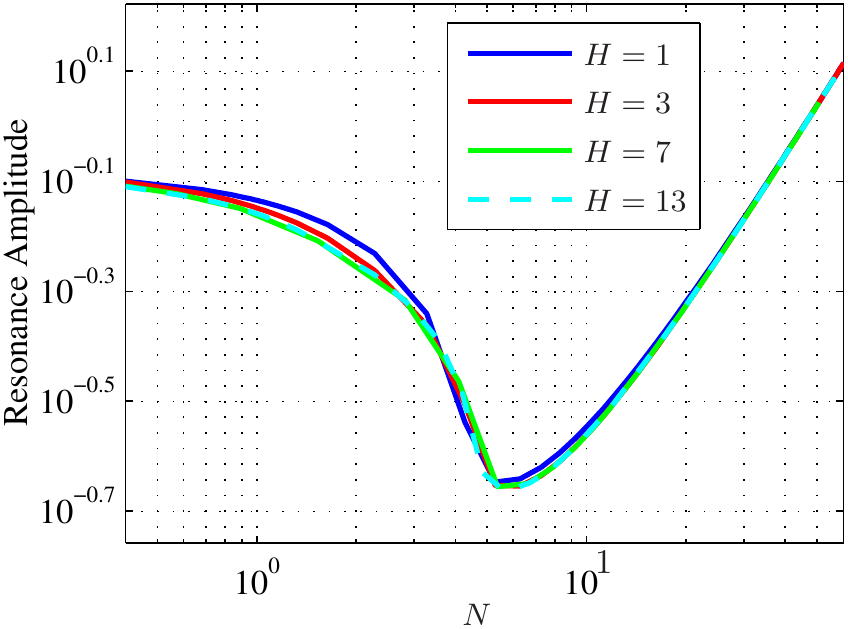}{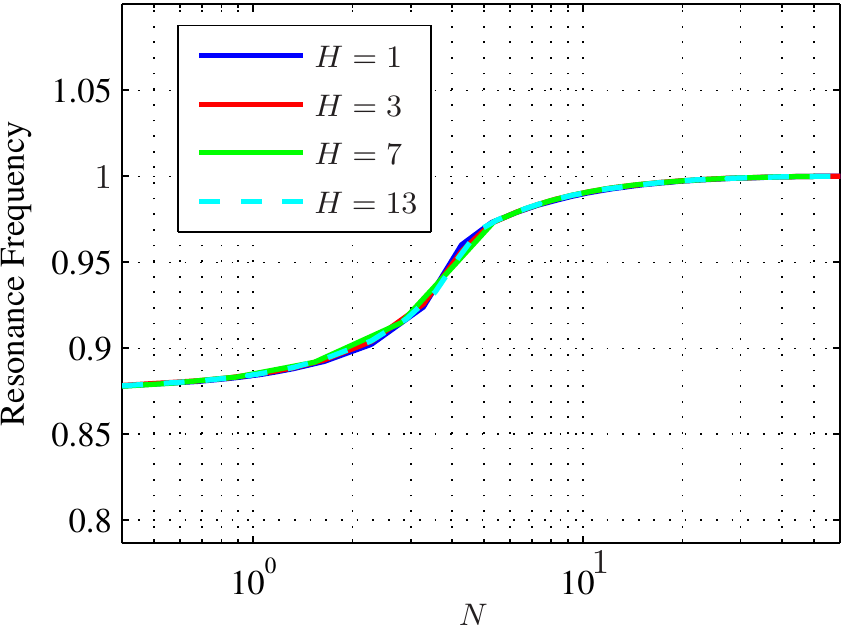}{}{}{.45}{.45}{Resonance properties for
different orders of the harmonic balance approach (~(a) resonance
amplitude, (b) resonance frequency~)}
\\
A harmonic force excitation was imposed at the center of the free
end in a frequency range near the second bending eigenfrequency. The
forced response for varied normal preload \g{N=-g/\kn} is depicted
in \fref{fig12a} along with the backbone curve. The results are
generally similar to the ones presented in \ssref{twodoffric}.
Again, it can be ascertained that there exists an optimum normal
preload that minimizes the resonance amplitude. For smaller preload
values, the contact node may lift off during one period of
oscillation. This causes a softening effect, leading to overhanging
branches in the forced response characteristic. In \fref{fig12b},
some periodic orbits corresponding to points on the backbone curve
are illustrated in a three-dimensional section through the phase
space. According to expectations, a multiharmonic character can be
ascertained from the response. In particular, the static component
of the displacement is varying with the vibration
amplitude.\\
In \frefs{fig13b}-\frefo{fig13a}, resonance amplitude and frequency
are depicted as a direct function of the normal preload \g{N}.
Several harmonics are required to achieve asymptotic convergence of
the resonance properties. This particularly holds for smaller values
of \g{N} \ie when the oscillation of the normal load and partial
separation gain influence on the dynamics of the system.

\section{Conclusions\label{sec:conclusions}}
A method was proposed that allows for an analytical formulation of
the high-order Harmonic Balance Method for the dynamic analysis of
systems with distinct states. The method can be applied to generic
nonlinearities that can be described by piecewise polynomial
functions and state transition conditions.\\
The methodology not only facilitates the computation of the periodic
solution but also provides accurate sensitivity data of the solution
to arbitrary system parameters that can be used \eg for design
studies. It was shown that the approach can be superior to the
conventional Alternating-Frequency-Time scheme in terms of accuracy
and computational efficiency, in particular if the sensitivities of
the transition time instants between the states are of
interest.\\
The method was applied to several structural dynamical systems with
conservative and dissipative nonlinearities in externally excited
and autonomous configurations. Generally good performance and
robustness of the numerical method were observed.\\
Possible future work includes a comparison of the method to the
Harmonic Balance formulation of the Asymptotic Numerical Method, as
introduced in \zo{coch2009}, and the application to other
engineering fields such as electrical switching networks.

\begin{appendix}
\section{Definite integral of a truncated Fourier series\label{asec:fourier_integral}}
A truncated Fourier series $a(\tau)$ is considered,
\e{a(\tau) = \suml{m=-H}{H}A_m\ee^{\ii m\tau}\fp}{a1}
Substituting this definition into \eref{fourier_integral} yields
\e{\intl{\tau^-}{\tau^+} a(\tau)\ee^{-\ii n\tau}\dd\tau =
\intl{\tau^-}{\tau^+}\suml{m=-H}{H}A_m\ee^{\ii m\tau}\ee^{-\ii
n\tau}\dd\tau = \suml{m=-H}{H}A_m\intl{\tau^-}{\tau^+}\ee^{\ii
(m-n)\tau}\dd\tau\fp}{a2}
The indefinite integral of the integral in the last part of
\eref{a2} can be expressed as
\e{\int\ee^{\ii (m-n)\tau}\dd\tau = \begin{cases}\tau &
m=n\\
\frac{\ee^{\ii(m-n)\tau}}{\ii(m-n)} & m\neq n\end{cases}\fp}{a3}
The case $m=n$ thus has to be treated separately. For convenience,
the sum in \eref{a2} is therefore split up. Substituting \eref{a3}
into \eref{a2} finally gives
\eal{\intl{\tau^-}{\tau^+} a(\tau)\ee^{-\ii n\tau}\dd\tau
&=&\left[A_n\tau + \suml{m=-H,m\neq
n}{H}A_m\frac{\ee^{\ii(m-n)\tau}}{\ii(m-n)}\right]_{\tau^-}^{\tau^+}\\
&=& \left(\tau^+-\tau^-\right) A_n + \suml{m=-H,m\neq
n}{H}\frac{\ee^{\ii\left(m-n\right)\tau^+}-
\ee^{\ii\left(m-n\right)\tau^-}}{\ii\left(m-n\right)}A_m\fp}{a4}
\end{appendix}


%% References


%\bibliographystyle{elsarticle-num}
%\bibliography{literature_mk}

\end{document}